\documentclass[11pt]{article}
\usepackage{amsmath,amsfonts,amssymb,latexsym,amsbsy, bbm, theorem,enumerate,color}
\usepackage{graphicx, graphics}
\usepackage{float,color,fancybox,shapepar,setspace,hyperref}
\usepackage[affil-it]{authblk}
\usepackage{caption}
\usepackage{tikz}
\usetikzlibrary{arrows}
\usepackage[normalem]{ulem}
\usetikzlibrary{arrows,automata}
\usepackage[latin1]{inputenc}
\usepackage{verbatim}

\usepackage{subfigure}
\usepackage{psfrag}
\usepackage{epstopdf}

\newtheorem{thm}{Theorem}[section]

\newtheorem{cor}[thm]{Corollary}%[section]
\newtheorem{prop}[thm]{Proposition}%[section]
\theorembodyfont{\rmfamily}
\begin{document}

\title{On edge-weighted mean eccentricity of graphs}

\author{Peter Johnson and Fadekemi Janet Osaye\thanks{Some of the first result presented in this paper forms part of the  author's PhD thesis at the University of Johannesburg, South Africa. Financial support from 2018 Kovaleskaia Research Grants for Female Mathematicians is gratefully acknowledged.}}

%Joshua Carlson$^2$, Andrew Owens$^3$,  K. E. Perry$^4$\\ Inne Singgih$^5$,  
%Zi-Xia Song$^{6,}$\thanks{Partially supported by the National Science  Foundation under Grant No. DMS-1854903.}\,,  Fangfang Zhang$^{7,}$\thanks{Corresponding author. This work was done in part while   Fangfang Zhang visited the University of Central Florida as a visiting student.   The visit was 
%supported by the Chinese Scholarship Council.    E-mail addresses: dlyang120@163.com (D. Yang);     jc31@williams.edu (J. Carlson);   adowens@widener.edu (A.  Owens);  kperry@soka.edu (K. E. Perry);  inne.singgih@uc.edu (I. Singgih); Zixia.Song@ucf.edu (Z-X. Song);  Fangfangzh@smail.nju.edu.cn (F. Zhang);  zhangx42@myumanitoba.ca (X. Zhang).}\,,  Xiaohong Zhang$^8$
  \affil{
  { \small {Department  of Mathematics and Statistics, Auburn University, AL  36830, AL}}\\
   }

% \date{ }

\maketitle
\begin{abstract}
Let $G$ be a connected edge-weighted graph of order $n$ and size $m$. Let $w:E(G)\rightarrow \mathbb{R}^{\geq 0}$ be the weighting function. We assume that $w$ is normalised, that is, $\sum_{e\in E(G)} w(e)=m$. The weighted distance $d_w(u,v)$ between any two vertices $u$ and $v$ is the least weight between them and the eccentricity $e_w(v)$ of a vertex $v$ is the weighted distance from $v$ to a vertex farthest from it in $G$. The mean(average) eccentricity of $G$, $avec(G,w)$, is the (weighted) mean of all eccentricities in $G$.  We obtain upper and lower bounds on $avec(G,w)$ in terms of $n$, $m$ or edge-connectivity $\lambda$ for two cases: $G$ is a tree and $G$ is connected but not a tree. In addition, we obtain the Nordhaus-Gaddum-type results for edge-weighted average eccentricity.
\end{abstract}

\baselineskip 18pt
\section{Introduction}

A graph is $G$ is said to be \textit{edge-weighted} if there is a weight function $w: E(G)\rightarrow \mathbb{R}^{\geq 0}$ which assigns to every edge $e\in E(G)$ a nonnegative function $w(e)$ called the \textit{weight} of $e$. An edge-weighted graph may be thought of as a special case of multigraphs, $M$ with no loops, such that the weight on each edge of $G$ can be interpreted as the multiplicity of the edge between two vertices $u$ and $v$ in $E(M)$. The weighted distance $d_w(u,v)$ is the least weight between any two vertices $u$ and $v$ and the eccentricity $e_w(v)$ of the vertex $v$ is the weighted distance between $v$ and a vertex farthest from $v$ in $G$.
\\
\\
If a graph models a routing network, for example, with vertices representing routers, edges representing links between the routers, and weights representing network packets on each link that can be forwarded from one router to the other, then the eccentricity of a vertex could be an indicator for a router's worst-case packet quantity. A minimum average eccentricity indicates the routers carry relatively low packets compare to the total packets for distribution, whereas a high average eccentricity implies that the routers carry high amounts of packets relative to the total packets. However, a redistribution of the network packets may increase or decrease an average eccentricity. Often, trees are the cheapest network to consider. We investigate how the distribution of weights affect the parameter when one considers trees of given number of end vertices. In addition, it is desirable to monitor a network's average eccentricity and that of its alternative links for optimal result. That is, what is the least or maximum difference that can occur when one compares the average eccentricity for a graph and that its complement? 

The average eccentricity of unweighted graphs, that is, graphs whose edges have all unit length, was first introduced by Buckley and Harary in \cite{BucHar} under the name \textit{eccentric mean} but only attracted much attention after its first systematic study by Dankelmann et.al. in \cite{DanGodSwa2004}, who amongst other results, determined the maximum average eccentricity of an unweighted connected graph of given order.
\begin{prop}\cite{DanGodSwa2004} \label{prp}
	If $G$ is a connected unweighted graph of order $n$, then
	\begin{equation*}
	{\rm avec}(G)\leq \frac{3n}{4}-\frac{1}{2},
	\end{equation*}
	with equality if and only if $G$ is a path of even order.
\end{prop}
In the same paper, the Nordhaus-Gaddam result for the unweighted average eccentricity was obtained.
\begin{prop}\cite{DanGodSwa2004}
	Assume both $G$ and $\overline{G}$ are connected and $n\geq 5$. Then the following bounds are best possible. 
	$$4\leq {\rm avec}(G) +{\rm avec}(\overline{G}) \leq \frac{3n}{4}+\frac{3}{2},$$ $$4\leq {\rm avec}(G){\rm avec}(\overline{G})\leq \frac{3n}{2}-1.$$ 
\end{prop}
Since then, the average eccentricity has been studied by several other authors (see\cite{AliDanMorMukSwaVet2018},\cite{AlexDan2020}, \cite{DanOsa-manu},\cite{DanOsa2019},\cite{DanMukOsaRod},\cite{TZ} for references). One of the most notable result is the lower bound obtained by Tang and Zhou.
\begin{prop}\cite{TZ}\label{prp1} 
	Let $G$ be a unicyclic graph of order $n\geq 5$. Then
	\[{\rm avec}(G) \geq 2-\frac{1}{n}.  \]
	Equality holds if and only if $G$ is formed by adding one edge to the $n$-vertex star.
\end{prop}
 Results on some distance measures of an edge-weighted graph are known in the literature. However, not much is known about the average eccentricty of graphs when the weight on each edge exceeds one. Broere et. al. \cite{BroDanDor} studied the average distance of edge-weighted graphs when some properties of the graph are known. It was further shown in \cite{DjeKou} that if the collection of the weights to be assigned is fixed, then the problem of maximizing the average distance is NP-complete. They also showed that for $\lambda$-edge-connected multigraphs, the mean distance is bounded above by $ \frac{2m}{3\lambda}$. In this paper, we establish similar bounds on the average eccentricity of an edge-weighted graph if the structure of the graph is known but not the actual weight function. We first obtain a lower bound for trees given that $w$ is a normalised weight function and then for general graphs that are not trees. Finally, we establish the Nordhaus-Gaddum results for edge-weighted average eccentricity. 
 
 \section{Definitions and notations}
 Let $G=(V,E)$ be a connected  graph with $m$ edges and order $n$; and let $w:E(G) \rightarrow \mathbb{R}^{\geq 0}$ be a weight function. For $u,v\in V(G)$, we denote by $d_w(u,v)$, the minimum length of a path between $u$ and $v$ according to the valuation $w$. The eccentricity of a vertex $v$ in $G$ is the largest weighted distance from $v$ to any other vertex in $G$. Let $EX(G,w)$ be the sum of all weighted eccentricities in $G$. The average eccentricity of a weighted graph $G$, denoted by $avec(G,w)$, is defined to be the average of all weighted eccentricities in $G$, that is \[avec(G,w)=\frac{1}{n} \sum_{v\in V(G)} e_w(v).  \]
 A \textit{normalised} weight function is a nonnegative weight function whose average weight of edges equals 1, that is;
 \[\sum_{e\in E(G)}w(e) =m,\] where $m$ is the size of $G$ and $0\leq w(e)\leq m$. With this restriction, $avec(G,w)$ becomes a generalisation of the unweighted average eccentricity. On account of this, every nonnegative weight function can be obtained from a normalised weight function by multiplying with a suitable real number. If $w(e)=1$ for each $e \in E(G)$, we write $avec(G,1)$. Unless otherwise stated, every weight function used in this paper is normalised.
 
 We define the \textit{minimum average eccentricity} of a graph $G$ by $$avec_{\min}(G)= \inf_w avec(G,w),$$ and the \textit{maximum average eccentricity} of $G$ by $$avec_{\max}(G)= \sup_w avec(G,w),$$ where the infimum and supremum are taken over all normalised weight functions $w$ of $G$. 
 A normalised weight function for which $avec_{\max}(G) = avec(G,w)$ is called an \textit{optimal} weight function of $G$.  Given an edge-weighted graph $G$, bounds on $avec(G,w)$ are uninteresting if we have no further information or restriction on $w$. This is because the multiplication of the weight function $w$ by a constant factor $c$ can yield an arbitrarily large or small value for $avec(G,cw)$; depending on the value of $c$.
 
 If $H$ is a subgraph of $G$, we write $w(H)$ for $\sum_{e\in E(H)}w(e)$. The degree of a vertex $v$ is the number of vertices adjacent to $v$ in $G$. Let $T$ be a tree. A \textit{leaf} or \textit{end vertex} is a vertex of degree 1. An \textit{end edge} in $T$ is an edge incident with a leaf of $T$ and every edge that is not an end edge of $T$ is an \textit{internal edge} of $T$. An edge $e\in E(G)$ is called an \textit{edge-cut} if $G-e$ disconnects $G$. A subset $F$ of $E$ is a \textit{cut set} if deleting all edges in $F$ from $G$ disconnects $G$. The number of edges in a smallest cut set of $G$ is the \textit{edge-connectivity} of $G$, denoted as  $\lambda$. A graph $\overline{G}$ is the \textit{complement} of $G$ on same vertices in which two distinct vertices of $\overline{G}$ are adjacent if and only if they are not adjacent in $G$. A graph is said to be \textit{self-complementary} if it is isomorphic to its complement.
 
 We define a \textit{\textbf{single broom}} $B(3,n-3)$ or simply $B_n$, as the tree obtained from $P_3$ by appending a set of $n-3$ end vertices labelled $A= v_1, v_2, \ldots, v_{n-3}$ to an end vertex of the $P_3$.
 
 \section{Main Results}
 \subsection{Edge-weighted trees}
 Let $t$ be the number of leaves in a tree $T$. The following theorem shows that for $T$, a weight function minimises the total eccentricity if we assign a zero weight to all internal edges of $T$ and weight $\frac{n-1}{t}$ to the end edges of $T$. 
 \begin{thm}\label{thm4-1}
 	Let $T$ be a tree with $n$ vertices, $n\geq 3$, and $t$ end vertices. If 
 	$w$ is a normalised weight function on the edges of $T$, then 
 	\[ {avec}(T,w)  \geq  \frac{(n-1)(n+t)}{nt}. \]
 	Equality holds if and only if every end edge has weight $\frac{n-1}{t}$,
 	and all the other edges have weight $0$.
 \end{thm}
\textbf{Proof.} Let $u_1, u_2, \ldots, u_t$ be the end vertices of $T$. For $i=1,2,\ldots,t$ let $e_i$ 
be the end edge incident with $u_i$. \\
Let $w$ be an arbitrary normalised weight function of $T$. We first show that 
for every internal vertex $v$ of $T$
\begin{equation}  \label{eq:internal-vertex-ecc}
e_w(v) \geq \frac{n-1}{t} \quad 
\textrm{for all $v\in V(T)$.} 
\end{equation}
Since each edge of $T$ is on a $(v,u_i)$-path for some $i\in \{1,2,\ldots,t\}$, 
we have 
\[ \sum_{i=1}^t d(v,u_i) 
\geq \sum_{e \in E(T)} w(e)  
= n-1. \]
Hence there exists $i\in \{1,2,\ldots,t\}$ such that 
\[  \frac{n-1}{t} \leq d(v,u_i) \leq e_w(v), \]
and so \eqref{eq:internal-vertex-ecc} follows. We now bound the sum of the eccentricities of the end vertices.                  
The eccentricity of an end vertex $u_i$ is clearly at least 
the average of the distances from $u_i$ to the other end 
vertices $u_j$, $j \neq i$, that is,
\[   e_w(u_i) \geq \frac{1}{t-1} \sum_{j: j\neq i} d(u_i, u_j). \]
Summing this over all $i\in \{1,2,\ldots,t\}$ yields
\begin{equation}  \label{eq:weighted-lower-bound-1} 
\sum_{i=1}^t e_w(u_i) \geq 
\frac{1}{t-1} \sum_{i=1}^t \sum_{j: j\neq i} d(u_i, u_j). 
\end{equation}
For every edge $e$ of $T$ let $c(e)=2n_1 n_2$, where $n_1$ and $n_2$ 
are the number of end vertices of $T$ in the two components of $T-e$. 
Since $c(e)$ counts the number of ordered pairs $(u_i, u_j)$ 
for which the $(u_i, u_j)$-path in $T$ contains $e$, we have 
\begin{equation}  \label{eq:weighted-lower-bound-2}
\sum_{i=1}^t \sum_{j: j\neq i} d(u_i, u_j) 
= \sum_{e\ \in E(T)} c(e) w(e). 
\end{equation}  
If $e$ is an end edge, then $c(e) =2(t-1)$, and if $e$ is not an end edge
then $c(e) \geq 4(t-2)$. So $c(e) \geq 2(t-1)$ with equality if 
and only if $e$ is an end edge. 
Since $\sum_{e\in E(T)} w(e) = n-1$, it follows that  
\begin{equation}   \label{eq:weighted-lower-bound-3}
\sum_{e \in E(T)} c(e) w(e) \geq \sum_{e \in E(T)} 2(t-1)w(e)
= 2(t-1)(n-1). 
\end{equation}  
Equality in \eqref{eq:weighted-lower-bound-3} implies that 
only the end edges of $T$ have positive weight. 
From \eqref{eq:weighted-lower-bound-1}, \eqref{eq:weighted-lower-bound-2}
and \eqref{eq:weighted-lower-bound-3} we get that
\begin{equation} \label{eq:end-vertices-ecc}
\sum_{i=1}^t e_w(u_i) \geq  2(n-1).
\end{equation}
Now from \eqref{eq:internal-vertex-ecc} and \eqref{eq:end-vertices-ecc},
we get that
\[ EX(T,w) \geq (n-t)\frac{n-1}{t} + 2(n-1) = \frac{(n-1)(n+t)}{t}, \]
and division by $n$ yields the bound in the theorem.   \\
\\
To prove the second part of the theorem assume that the bound holds with equality. 
Then we have
equality in \eqref{eq:end-vertices-ecc} and in \eqref{eq:internal-vertex-ecc}. 
Equality in \eqref{eq:end-vertices-ecc} implies equality in
\eqref{eq:weighted-lower-bound-3}, and thus that only end edges
have positive weight. Equality in \eqref{eq:internal-vertex-ecc} implies
that no edge has weight strictly greater than $\frac{n-1}{t}$, which
in turn implies that all end edges have weight $\frac{n-1}{t}$. 
Finally, assume that $w$ is the weight function defined by 
$w(e)=\frac{n-1}{t}$ for every end edge, and $w(e)=0$ if $e$ is
not an end edge.  Then $e(v)= \frac{2(n-1)}{t}$ if $v$ is an end vertex
of $T$, and $e(v)=\frac{n-1}{t}$ if $v$ is an internal vertex of $T$, so
${ avec}(T,w)= \frac{(n-1)(n+t)}{nt}$.

\hfill $\Box$ 
\begin{thm}\label{props4-2}
	Let $T$ be a tree of order $n \geq 2$ and $w$ a normalised weight
	function on $E(T)$. Then 
	\[ { avec}(T,w) \leq n-1. \]
	Equality holds if and only if there exists an edge 
	$e$ of $T$ with $w(e)=n-1$.
\end{thm}
 \textbf{Proof.} 
 Let $w$ be an arbitrary normalised weight function on $T$. 
 The total weight of all edges of $T$ is $n-1$. Hence no
 distance between two vertices in $T$ is greater than $n-1$, and thus
 \[ e_w(v) \leq n-1 \quad \textrm{for all $v \in V(T)$.} \] 
 Summing over all $v\in V(T)$ and dividing by $n$ yields 
 \[ { avec}(T,w) \leq n-1. \]
 To prove the second part of the theorem, assume that $w$ is 
 a normalised weight function on $E(T)$ for which ${avec}(T,w) = n-1$.
 Then $e_w(v) =n-1$ for all $v\in V(T)$. 
 We show that there exist an edge of $T$ of weight $n-1$.
 Suppose not. Then there exist two edges $e_1, e_2$ of positive
 weight in $T$. Let $e_1=uv$, and assume without loss of generality
 that $u$ is closer to $e_2$ than $v$. 
 Let $u'$ be an eccentric vertex of $u$ and let $P$ be the $(u,u')$-path
 in $T$. Then $P$ cannot contain both, $e_1$ and $e_2$. If it does not
 contain $e_1$, then 
 $$d(u,u') = \sum_{e \in E(P)} w(e) \leq n-1-w(e_1) < n-1.$$ 
 Similarly, if $P$ does not contain $e_2$, then 
 $$d(u,u') = \sum_{e \in E(P)} w(e) \leq n-1-w(e_2) < n-1.$$ 
 In both cases we obtain the contradiction $e_w(u) < n-1$. 
 Hence there exists an edge $e$ with $w(e)=n-1$. 
 On the other hand, if $w$ is a weight function that assigns weight
 $n-1$ to one edge and $0$ to all other edges, then clearly every vertex has eccentricity $n-1$, and so ${\rm avec}(T,w)=n-1$. 
 \hfill $\Box$
\\
\\ 
The following corollary is a consequence of Theorems \ref{thm4-1} and  \ref{props4-2}.

\begin{cor}\label{cor4-1}
	Let $T$ be a tree of order $n\geq 3$. Then
	\[2-\frac{1}{n} \leq avec_{\min}(T) \leq \frac{n+1}{2}-\frac{1}{n},\]  and \[avec_{\max}(T) = n-1.\]
\end{cor}
\textbf{Proof.} 
Since the tree $T$ with the minimum number of end vertices is the path $P_n$ with $t=2$ and $T$ with maximum number of end vertices is the star $S_n$ with $t=n-1$,
it follows by Theorem \ref{thm4-1} that the lower bound for $avec_{min}(T)$ is attained when $t=2$ and the upper bound is attained when $t=n-1$. By substituting for $t=2$ and $t=n-1$, respectively, the first part of the Corollary \ref{cor4-1} follows.
\\
To prove the second part of Corollary \ref{cor4-1}, we maximize $EX(T,w)$ by moving all weight onto any edge $e''$ by Theorem \ref{props4-2}. 
Then $e_w(v) =n-1$, for every vertex $v$ in $T$ and thus $EX_{max}(T) =n(n-1)$ and so $avec_{max}(T) = n-1$.
\hfill $\Box$
\\
\\
We remark that both inequalities of the first part of Corollary \ref{cor4-1} are best possible, since equality is attained on the left hand side by the star $K_{1,n-1}$ and on the right hand side by the path $P_n$.
%In Theorem \ref{thm4-1} above, we proved that the lower bound also applies for weighted trees, while the upper bound has to be increased by a value of $\frac{n}{4}-\frac{1}{2}$.  
Corollary \ref{cor4-1} generalises the lower bound on the average eccentricity of unweighted trees in Proposition \ref{prp1}. However, it does not generalise the upper bound in Proposition \ref{prp} but differs by $\frac{n}{4}-\frac{1}{2}$.

\subsection{Edge-weighted graphs}
We now turn our attention to connected graphs that are not trees. 
\begin{thm}\label{propp}
	If $G$ is a connected graph but not a tree, with $m$ edges and edge-connectivity $\lambda$ then $$0\leq avec(G,w)\leq \frac{m}{\lambda}.$$ 
\end{thm}
\textbf{Proof.}  
We first prove the lower bound. Clearly we have $avec_{min}(G) \geq 0$ since we consider only nonnegative weight functions for which every distance is nonnegative. 
It suffices to construct a weight function $w$ with $EX(G,w)=0$. If $G$ is not a tree, then $m\geq n$. Let $T$ be a spanning tree of $G$. Let $w$ be a weight function that assigns weight $0$ to all edges in $T$, and the remaining edges receive an arbitrary nonnegative weight so that 
the total weight is $m(G)$. Since $T$ contains a path of weight $0$ between any two vertices of $G$, we have $e_w(v)=0$ for all $v\in V(G)$, and thus $EX(G,w)=0$. Hence the minimum possible average eccentricity is zero.
\\
\\
To prove the upper bound of the theorem, let $w$ be an arbitrary normalised weight function on $E(G)$. By Menger's Theorem there are $\lambda$ edge-disjoint $u-v$ paths in $G$ for all $u,v \in V(G)$. Among these there exists a path of total weight at most $\frac{m}{\lambda}$. Hence $d_w(u,v) \leq \frac{m}{\lambda(G)}$ for all $u,v\in V(G)$, which implies that $e_w(v) \leq \frac{m}{\lambda}$ for all $v\in V(G)$. 
Hence $EX_w(G) \leq \frac{mn}{\lambda}$ and thus $avec_w(G) \leq \frac{m}{\lambda}$. Since $w$ was arbitrary, we conclude that
\[ avec(G,w) \leq \frac{m}{\lambda}. \]
To show that $avec_{max}(G) = \frac{m}{\lambda}$ it suffices to give a weight function $w$ for which $avec(G) = \frac{m}{\lambda}$. 
Let $S$ be a minimum edge-cut of $G$. Define the weight function $w$ by 
\[ w(e) = \left\{\begin{array}{cc}
\frac{m}{\lambda} & \textrm{if $e \in S$,} \\
0  & \textrm{if $e \notin S$} 
\end{array} \right. \]
Clearly, $w$ is normalised. Also, the distance between any two vertices in the same component of $G-S$ is $0$, while the distance between two vertices in different components of $G-S$ is $\frac{m}{\lambda}$. Hence $e_w(v)=\frac{m}{\lambda}$ for all $v\in V(G)$ and thus $avec_w(G) = \frac{m}{\lambda}$, as desired. 
\hfill $\Box$
\\
\\
{\textbf{Remark}:}\,\, We have seen that the average eccentricity of a weighted graph can be much larger than that of the underlying unweighted graph. The above theorem is a generalisation of the corresponding result for trees, whose weighted maximum average eccentricity is always $n-1$. It is important to note that the result for trees, which says that the weight  functions realising the maximal average eccentricity are those in which all the weight is concentrated in one edge, does not generalise to general graphs. For example the maximum average eccentricity of an even cycle can also be realised, among others, by the all-1 function.

\subsection{Nordhaus-Gaddum Type}
As a consequence of the results obtained for ${\rm avec}(G,w)$, the Nordhaus-Gaddum results for edge-weighted average eccentricity are easy to obtain. We consider cases for trees and connected graphs that are not trees. First, we prove the lower and upper bounds on the sum $avec(G,w) + avec(\overline{G},w)$ when $G$ is a tree and characterise the associated extremal conditions.
\begin{cor}\label{cor-1}
	Let $G$ be an $n$-vertex tree with $n\geq 4$. \\
	(a)\,\ 	If $\overline{G}$ is also a tree, then $G=\overline{G} =P_4$. Thus,  $$\frac{9}{2} \leq avec(G,w)+ avec(\overline{G},w) \leq 6,\quad \text{and} \quad \frac{81}{16}\leq avec(G,w)avec(\overline{G},w)\leq 9.$$
	(b)\,\, If $\overline{G}$ is connected but not a tree with $n\geq 5$, size $\overline{m}$ and edge-connectivity $\overline{\lambda}$, then \[\frac{2(n-1)^2}{n(n-2)} \leq avec(G,w) +avec(\overline{G},w)\leq \frac{(n-1)[n+2(\overline{\lambda}-1)]}{2\overline{\lambda}}\] {and}\[ 0\leq avec(G,w)avec(\overline{G},w)\leq \frac{\overline{m}(n-1)}{\overline{\lambda}}.\] The lower bound is attained if $G$ has $n-2$ end vertices and the upper bound is attained  if $G$ is a single broom.               
\end{cor}
{\bf Proof:}\\
(a) Generally, it is trivial that since the sum of the sizes of $G$ and $\overline{G}$ on a given order equals the size of the complete graph of the same order, and since $G$ is a tree, we have $|E(\overline{G})|=\binom{n}{2}-(n-1)$ which when simplified yields $\frac{(n-1)(n-2)}{2}$. Since $\overline{G}$ is also a tree with $n-1$ edges, solving for $n$ in $\frac{(n-1)(n-2)}{2} = n-1$ yields that $n=4$ and $n=1$. It follows that amongst all trees of order $4$, the path $P_4$ is the only self-complementary tree. By Theorem \ref{thm4-1}, the lower bound is attained with $n=4$ and $t=2$ and by Theorem \ref{props4-2}, the upper bound is attained since $avec_{max}(P_4) =3$.\\
\\
(b) We first prove the lower bound. Since $G$ is a tree and $\overline{G}$ is connected but not a tree, $|E(\overline{G})|= \frac{(n-1)(n-2)}{2}$ by the same argument as in the proof of (a). 
%Thus by Theorem \ref{thm4-1}, $$ avec(G) \geq \frac{(n-1)(n+t)}{nt},$$ where $t$ is the number of end vertices of $G$. 
It is easy to see that $2\leq t\leq n-2$ since letting $t=n-1$ results in a $\overline{G}$ that is disconnected. Denote $\frac{(n-1)(n+t)}{nt}$ by $f(n,t)$. A straight forward minimisation shows that $f(n,t)$ is minimised when $t=n-2$. 
Hence, 
\begin{equation*}
\begin{aligned}
avec(G,w) &\geq \frac{(n-1)(2n-2)}{n(n-2)}\\
&= \frac{2(n-1)^2}{n(n-2)}. 
\end{aligned}
\end{equation*}
Thus, since $avec(\overline{G},w) \geq 0$ by Theorem \ref{propp}, it follows that 
\[avec(G,w) +avec(\overline{G},w)\geq \frac{2(n-1)^2}{n(n-2)}. \]
\\
\\
To prove the upper bound, let $\overline{m}$ and $\overline{\lambda}$ be the size and edge-connectivity of $\overline{G}$, respectively. We use the fact that $avec_{max}(G) = n-1$ from Corollary \ref{cor4-1} and $avec_{max}(\overline{G}) = \frac{\overline{m}}{\overline{\lambda}}$ from Theorem \ref{propp}. Note that $\overline{m} = \binom{n}{2}-(n-1) = \frac{(n-1)(n-2)}{2}$. It then follows that
\begin{eqnarray*}
	avec(G,w) +avec(\overline{G},w) &\leq &	n-1 + \frac{\overline{m}}{\overline{\lambda}} \\
	&=& (n-1) + \frac{(n-1)(n-2)}{2\overline{\lambda}} \\
	%&\leq & (n-1) + \frac{(n-1)(n-2)}{2}\\
	&=& \frac{(n-1)[n+2(\overline{\lambda}-1)]}{2\overline{\lambda}}.
\end{eqnarray*}
\\
\\
To prove that the bounds are sharp, let $P_3$ be the path labelled with end vertices $a$ and $c$, and center vertex $b$. Let $B_n$ be the single broom obtained from $P_3$ by appending to $c$, a set of $n-3$ end vertices labelled $A= v_1, v_2, \ldots, v_{n-3}$. Then $B_n$ has $n-2$ end edges and an internal edge $bc$. Let $\overline{B}_n$ be the complement of $B_n$ with size $\overline{m}$. Since $c$ is a vertex of degree $n-2$ in $B_n$,  $deg_{\overline{B}_n}(c)=1$ so that the edge $ac$ is a bridge in  $\overline{B}_n$ and every other vertices in $\overline{B}_n$ is connected to at least $n-3$ vertices. Hence $\overline{B}_n$  is connected and contains a complete graph $K_{n-3}$ as an induced subgraph. 
\\
\\
To minimise $avec(B_n,w)$ we assign weights $\frac{n-1}{n-2}$ to each end edge and weight $zero$ to the edge $bc$ so that every end vertex has eccentricity $2(\frac{n-1}{n-2})$ and $e_w(b)=e_w(c)=\frac{n-1}{n-2}$. To minimise $avec(\overline{B},w)$, we assign weight zero to the edges of a spanning tree, $T$ of $\overline{B}$ while the  remaining edges not in $T$ receive arbitrary weight such that the total weight equals $\overline{m}$.
It follows that 
\[EX(B_n,w)+ EX(\overline{B}_n,w)\geq \frac{2(n-1)}{n-2}(n-2) + \frac{2(n-1)}{n-2}=\frac{2(n-1)^2}{n-2}.  \]
A division by $n$ yields the lower bound.\\
\\
To maximise $avec(B_n,w)$, we move all weights to any edge so that the eccentricity of every vertex becomes $n-1$. To maximise $avec(\overline{B}_n,w)$, we assign weight $\overline{m}$ to the edge $ac$ since $\overline{B}_n$ is 1-edge connected. Hence, $d_w(v,a) =\frac{(n-1)(n-2)}{2}$ for all $v\in V(G)-\{a\}$ implying that $e_w(v)=\frac{(n-1)(n-2)}{2}$ for all $v\in V(G)$. Thus $avec(\overline{B}_n,w) = \frac{(n-1)(n-2)}{2}$, and hence,
\begin{eqnarray*}
	avec(G,w) +avec(\overline{G},w)
	&=& (n-1) + \frac{(n-1)(n-2)}{2}= \frac{n(n-1)}{2}.
\end{eqnarray*}
This shows that equality is attained when $\overline{G}$ has a bridge.
\hfill $\Box$
\begin{figure}	
	\begin{center}
	\begin{tikzpicture}
	[scale=0.6,inner sep=0.7mm, % this is the node radius
	vertex/.style={circle,thick,draw}, % this defines the default style for the vertex class
	thickedge/.style={line width=3pt}] % this defines the default style for the thick edges
		\begin{scope}[>=]   
		\node[vertex] (a1) at (1,0) [fill=black] {};           
		\node[vertex] (a2) at (4,0) [fill=black] {};  
		\node[vertex] (a3) at (7,0) [fill=black] {};
		\node[vertex] (a4) at (10,1.3) [fill=black] {};
		\node[vertex] (a5) at (10,-1.3) [fill=black] {};
		\node[vertex] (a6) at (8.5,2.7) [fill=black] {};
		\node[vertex] (a7) at (8.5,-2.7) [fill=black] {};
		\draw[very thick] (a3)--(a4); 
		\draw[very thick] (a3)--(a5);
		\draw[ very thick] (a3)--(a6);
		\draw[very thick] (a3)--(a7);
		\draw[very thick] (a3)--(a2);
		\draw[very thick] (a3)--(a1);
%		%	\draw[very thick] (a1)--(a13);
%		%\draw[very thick] (a1)--(a14);
%		%\draw[ dashed] (a11)--(a12);
	\node[right] at (0.7,-0.4) {$a$};
	\node[right] at (3.7,-0.4) {$b$};
	\node[right] at (6.7,-0.4) {$c$};
	\node[right] at (9.7,0.9) {$v_2$};
	\node[right] at (9.7,-1.7) {$v_3$};
	\node[right] at (8.3,2.3) {$v_1$};
	\node[right] at (8.3,-3.1) {$v_4$};
		\end{scope}         
		\end{tikzpicture} 
	\end{center}
	\caption{Single Broom, $B(3,4)$ with $n=7$}
%	\label{fig:C4-free-sharpness-example}
\end{figure}
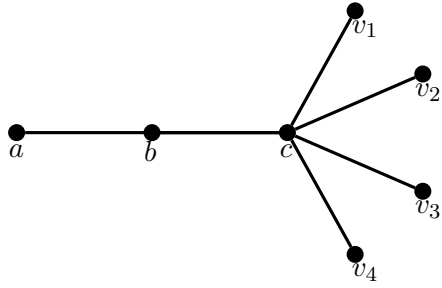

\begin{figure}	
	\begin{center}
		\begin{tikzpicture}
		[scale=0.6,inner sep=0.7mm, % this is the node radius
		vertex/.style={circle,thick,draw}, % this defines the default style for the vertex class
		thickedge/.style={line width=3pt}] % this defines the default style for the thick edges
		\begin{scope}[>=]   
	\node[vertex] (a1) at (1,0) [fill=black] {};           
	\node[vertex] (a2) at (4,0) [fill=black] {};  
	\node[vertex] (a3) at (1,-3) [fill=black] {};
		\node[vertex] (a4) at (4,-3) [fill=black] {};
		\node[vertex] (a5) at (-2,-1.5) [fill=black] {};
	\node[vertex] (a6) at (7,-1.5) [fill=black] {};
		\node[vertex] (a7) at (10,-1.5) [fill=black] {};
		\draw[very thick] (a1)--(a2); 
	    \draw[very thick] (a1)--(a3);
		\draw[ very thick] (a1)--(a4);
		
		\draw[very thick] (a2)--(a3); 
		\draw[very thick] (a2)--(a4);
		\draw[ very thick] (a3)--(a4);
		
		\draw[very thick] (a5)--(a1); 
		\draw[very thick] (a5)--(a2);
		\draw[ very thick] (a5)--(a3);
		\draw[very thick] (a5)--(a4);
		
		\draw[very thick] (a6)--(a1); 
		\draw[very thick] (a6)--(a2);
		\draw[ very thick] (a6)--(a3);
		\draw[very thick] (a6)--(a4);
		
		\draw[very thick] (a6)--(a7);
	\node[left] at (1.2,0.3) {$v_1$};
		\node[right] at (4.0,0.2) {$v_2$};
	\node[right of =a3] at (-0.65,-3.4) {$v_3$};
		\node[right] at (4.0,-3.2) {$v_4$};
		\node[left] at (-2,-1.7) {$b$};
		\node[right] at (7.1,-1.3) {$a$};
		\node[right] at (10.1,-1.3) {$c$};
		\end{scope}         
		\end{tikzpicture} 
	\end{center}
	\caption{Complement of the Single Broom, $\overline{B}(3,4)$ with $n=7$}
	%	\label{fig:C4-free-sharpness-example}
\end{figure}
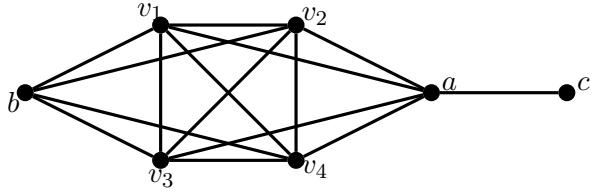

\begin{cor}
	  Let $G$ and $\overline{G}$ be connected non-tree graphs of edge-connectivities $\lambda$ and $\overline{\lambda}$. Then,
	  \[ 0\leq avec(G,w) + avec(\overline{G},w) \leq \frac{m}{\lambda}+\frac{\overline{m}}{\overline{\lambda}},   \] and 
	  \[0\leq avec(G,w)avec(\overline{G},w)\leq \frac{m\overline{m}}{\lambda \overline{\lambda}} .\]
\end{cor}
The proof of the corollary follows the same argument as Corollary \ref{cor-1}.
\hfill $\Box$
\subsection*{Acknowledgement}
The corresponding author most sincerely thank Dr. Peter Dankelmann for his immense support and mentorship during the course of her doctoral studies at the University of Johannesburg and for his helpful suggestions on this paper.  

%\noindent{\bf Acknowledgements}\bigskip 
%
%This work was completed in part at the 2019 Graduate Research Workshop in Combinatorics, which was supported in part by NSF grant \#1923238, NSA grant \#H98230-18-1-0017,  a generous award from the Combinatorics Foundation, and Simons Foundation Collaboration Grants \#426971 (to M. Ferrara), \#316262 (to S. Hartke) and \#315347 (to J. Martin).\bigskip

\end{document}